\documentclass[a4paper,11pt]{article}
\usepackage{amsmath,amscd,amsfonts,amssymb,epsf,latexsym}

\makeatletter

\long\def\@makefnt#1{\parindent 1em\noindent
            \hb@xt@1.8em{\hss\@textsuperscript{}}#1}
\long\def\@ftntext#1{\insert\footins{%
    \reset@font\footnotesize
    \interlinepenalty\interfootnotelinepenalty
    \splittopskip\footnotesep
    \splitmaxdepth \dp\strutbox \floatingpenalty \@MM
    \hsize\columnwidth \@parboxrestore
    \color@begingroup
      \@makefnt{%
        \rule\z@\footnotesep\ignorespaces#1\@finalstrut\strutbox}%
    \color@endgroup}}%
\def\subjclass#1{%
  \@ftntext{2000 {\itshape Mathematics Subject Classification.}\enspace #1.}}
\def\keywords#1{%
  \@ftntext{{\itshape Key words and phrases.}\enspace #1.}}
\makeatother

\def\moins{\raise 1pt\hbox{{$\scriptstyle -$}}}
\def\plus{\raise 1pt\hbox{{$\scriptstyle +$}} }

\newtheorem{theorem}{Theorem}
\newtheorem{proposition}[theorem]{Proposition}
\newtheorem{lemma}[theorem]{Lemma}
\newtheorem{corollary}[theorem]{Corollary}
\newtheorem{remark}[theorem]{Remark}

\def\proof{\noindent{\bf Proof.\ }}

\def\qed{~\hbox{$\Box$}}

\def\cJ{\mathop{\rm J}}

\def\codim{\mathop{\rm codim}}

\def\dim{\mathop{\rm dim}}
\def\Sym{\mathop{\rm Sym}}
\def\Hom{\mathop{\rm Hom}}

\def\corank{\mathop{\rm corank}}

\def\cT{{\mathcal T}}
\def\cJ{{\mathcal J}}

\begin{document}

\title{\bf Thom polynomials of invariant cones, Schur functions, and positivity}

\author{Piotr Pragacz\thanks{Research supported by a KBN grant.}\\
\small Institute of Mathematics of Polish Academy of Sciences\\
\small \'Sniadeckich 8, 00-956 Warszawa, Poland\\
\small P.Pragacz@impan.gov.pl
\and
Andrzej Weber\thanks{Research supported by the KBN grant 1 P03A 005 26.}\\
\small Department of Mathematics of Warsaw University\\
\small Banacha 2, 02-097 Warszawa, Poland\\
\small aweber@mimuw.edu.pl}

\date{(18.06.2007)}



\maketitle

\begin{abstract}
We generalize the notion of Thom polynomials from singularities of
maps between two complex manifolds to invariant cones in representations,
and collections of vector bundles. We prove that the generalized Thom polynomials,
expanded in the products of Schur functions of the bundles, have nonnegative
coefficients. For classical Thom polynomials associated with maps of complex
manifolds, this gives an extension of our former result for stable singularities
to nonnecessary stable ones. We also discuss some related aspects of Thom polynomials,
which makes the article expository to some extent.
\end{abstract}

\section{Introduction}\label{intro}

The present paper is both of the research and expository character. It concerns global
invariants for singularities. Our main new result here is Theorem \ref{Tm} (see also
Corollary \ref{polfun} and \ref{CpS}).

To start with, we recall that the global behavior of singularities
is governed by their {\it Thom polynomials} (cf. \cite{T},
\cite{AVGL}, \cite{Ka}, and \cite{Rim}). By a {\it singularity},
we shall mean in the paper a class of germs
$$
({\bf C}^{m},0) \to ({\bf C}^{n},0)\,,
$$
where $m,n \in {\bf N}$, which is closed under the
right-left equivalence (i.e. analytic
reparametrizations of the source and target).

Suppose that $f:M\to N$ is a map between complex manifolds,
where $\dim(M)=m$ and $\dim(N)=n$.
Let $V^{\eta}(f)$ be the cycle carried by the {\it closure} of the
set
\begin{equation}
\{x\in M : \hbox{the singularity of} \ f \ \hbox{at} \ x \
\hbox{is} \ \eta \}\,.
\end{equation}

We recall that the {\it Thom polynomial} ${\cal T}^{\eta}$ of a singularity
$\eta$ is a polynomial in the formal variables
$$
c_1, c_2,\ldots, c_m; \ c'_1, c'_2,\ldots, c'_n\,,
$$
such that after the substitution
\begin{equation}
c_i=c_i(TM), \ \ \  c'_j=c_j(f^*TN) \,,
\end{equation}
($i=1,\ldots, m$, $j=1,\ldots, n$) \
for a general map $f:M \to N$ between complex manifolds,
it evaluates the Poincar\'e dual\footnote{In the following, we shall often omit
the expression: ``the Poincar\'e dual of''.} of $[V^{\eta}(f)]$.
This is the content of the Thom theorem \cite{T}.
For a detailed discussion of the {\it existence} of Thom polynomials,
see, e.g., \cite{AVGL}. Thom polynomials associated with group actions
were studied by Kazarian in \cite{Ka}.

Recall that -- historically -- the first ``Thom polynomial''
appeared in the so-called ``Riemann-Hurwitz formula''.
Let $f:M\to N$ be a general  holomorphic map of compact Riemann surfaces.
This means that $f$ is a simple covering, that is, the critical
points are nondegenerate and at most one appears in each fiber.
Denoting by $e_x$ the {\it ramification index} of $f$ at $x\in M$
(i.e. the number of sheets of $f$ meeting at $x$), the
Riemann-Hurwitz formula asserts that
\begin{equation}
\sum_{x\in M} (e_x-1)=2g(M)-2-\deg(f)\bigl(2g(N)-2\bigr)\,.
\end{equation}
Denoting by $A_1$ the singularity of $z\mapsto z^2$ at $0$,
this is equivalent to saying that the fundamental class of the {\it ramification divisor}
of $f$\,,
$$
\sum_x (e_x-1)[x]=[V^{A_1}(f)]\,,
$$
where $x$ runs over the set of critical points of $f$, is given by the following
expression in the first Chern classes:
\begin{equation}
c_1(f^*TN)-c_1(TM)=c_1(f^*TN-TM)\,.
\end{equation}
In other words, the Riemann-Hurwitz formula says:
\begin{equation}
{\cal T}^{A_1}=c'_1-c_1\,.
\end{equation}
For a wider discussion of the Riemann-Hurwitz formula and early history of Thom polynomials,
we refer to Kleiman's survey article \cite{Kl}.
The Riemann-Hurwitz formula is true also in positive characteristic for finite separable morphisms
of algebraic curves (cf. \cite{H1}, Chap.~IV, Sect.~2).

\smallskip

Thom \cite{T} generalized the Riemann-Hurwitz formula to general
maps $f:M\to N$ of complex manifolds with $n-m>0$,
the singularity being always $A_1$:
\begin{equation}
[V^{A_1}(f)]=c_{n-m+1}(f^*TN-TM)=\sum_{i=0}^{n-m+1} S_{n-m+1-i}(TM^*) c_i(f^*TN)\,,
\end{equation}
where $S_j$ denotes the $j$th Segre class.

Though for the singularity $A_1$ the Thom polynomials are rather simple,
they start to be quite complicated even for simplest singularities
coming ``just after $A_1$'', say (cf. the tables in \cite{Rim}).
For example, the Thom polynomial
for the singularities $A_4$ is known only for small values of $k$
(cf. \cite{O}). Therefore, it is important to study the {\it structure} of Thom
polynomials. It appears that a good tool for this task is
provided by {\it Schur functions} \cite{FK}, \cite{P1}, \cite{P2},
\cite{P3}, \cite{P4}.
Let us quote some results related to Schur function expansions of Thom
polynomials of stable\footnote{By a {\it stable} singularity
we mean an equivalence class of stable germs $({\bf C}^{\bullet},0)
\to ({\bf C}^{\bullet+k},0)$, where $\bullet\in {\bf N}$, under the
equivalence generated by right-left equivalence and suspension
(by suspension of a germ $\kappa$ we mean its trivial unfolding:
$(x,v) \mapsto (\kappa(x),v)$). For a stable singularity, its Thom
polynomial is of the form $\sum_I \alpha_I S_I(TM^*-f^*TN^*)$,
where $S_I$ denotes a Schur function, cf. Sections 3 and 5.}
singularities.

First, it was shown in \cite{P2} that if a representative of a stable
singularity
$$
\eta:({\bf C}^{m},0) \to ({\bf C}^{n},0)
$$
is of Thom-Boardman type $\Sigma^i$, then all summands in the Schur function
expansion of ${\cal T}^{\eta}$ are indexed by partitions
containing\footnote{We say that
one partition {\it is contained} in another if this holds for their Young
diagrams.} the rectangle partition
$$
(n-m+i,\ldots,n-m+i) \ \ \ \ (i \ \hbox{times}).
$$
This is a consequence of the structure of the ${\cal P}$-{\it ideals of the singularities}
$\Sigma^i$, which were introduced and investigated in \cite{P-2}.
Second, in \cite{PW}, the authors proved that
for any partition $I$ the coefficient $\alpha_I$ in the Schur function expansion
of the Thom polynomial
$$
{\cal T}^{\eta}=\sum_I \alpha_I S_I(TM^*-f^*TN^*)
$$
is nonnegative. This result was conjectured before in \cite{FK}
and independently in \cite{P1}. It appears to be a consequence of
the Fulton-Lazarsfeld theory of numerical positivity of cones in
ample vector bundles \cite{FL} (cf. also \cite[\S 8]{Laz}),
combined with a functorial version of the bundles of jets,
appearing in the approach to Thom polynomials via classifying
spaces of singularities \cite{Ka}.

In the present paper, we shall prove a more general result that -- we believe --
will better explain a reason of the positivity in the above classical case,
as well as in many other situations. To this end, we extend the definition
of Thom polynomials from the singularities of maps $f: M\to N$ of complex
manifolds \cite{T} to the invariant cones in representations
of the product of general linear groups
$$
\prod_{i=1}^p GL_{n_i}\,.
$$
Such Thom polynomials are naturally defined on $p$-tuples of vector bundles
of ranks $n_i$. It is convenient to pass to {\it topological homotopy category},
where each $p$-tuple of bundles can be pulled back from the universal $p$-tuple
of bundles on the product of $p$ classifying spaces
$$
\prod_{i=1}^p BGL_{n_i}\,.
$$
Suppose that the functor associated with such a representation preserves
global generateness.
Our main result -- Theorem \ref{Tm} -- then asserts that the Thom polynomial
for a $p$-tuple of vector bundles $(E_1,E_2,\ldots,E_p)$, when expanded
in the basis
$$
\{S_{I_1}(E_1)\cdot S_{I_2}(E_2) \ \cdots \ S_{I_p}(E_p)\}
$$
of products of Schur functions applied to the successive bundles, has
{\it nonnegative} coefficients.
The key tool is {\it positivity of cone classes} for globally generated
vector bundles combined with the {\it Giambelli formula}.
For a polynomial representation of $\prod_{i=1}^p GL_{n_i}$ of positive degree,
we get, in addition, that the sum of the coefficients is {\it positive}
(cf. Corollary \ref{polfun}).

Theorem \ref{Tm}, in the classical situation of singularities of
maps $f:M\to N$ between complex manifolds, implies that for a given
singularity its Thom polynomial, when expanded in the basis
$$
S_I(TM^*)\cdot S_J(f^*TN)\,,
$$
has nonnegative coefficients (cf. Corollary \ref{CpS}).

We also note that Theorem \ref{Tm} implies the main result of our former paper
\cite{PW} for Thom polynomials of {\it stable} of singularities of maps between
complex manifolds, where, however, the Schur functions {\it in difference
of bundles} were used (cf. Theorem \ref{Tpos}).

\section{Thom polynomials of invariant cones}\label{Tpc}

In this section, we define ``generalized Thom polynomials''. Our construction
is modeled on that used to the construction of classical Thom polynomials with
the help of the ``classifying spaces of singularities'' (cf., e.g., \cite{Ka}).

\smallskip

Suppose that $(n_1,n_2,\ldots,n_p)\in {\bf N}^p$ and that $V$ is a representation of
\begin{equation}
G=\prod_{i=1}^p GL_{n_i}\,.
\end{equation}

The representation $V$ gives rise to a {\it functor} $\phi$ defined for
a collection of bundles on a variety $X$:
$$
E_1,E_2,\dots,E_p\mapsto \phi (E_1,E_2,\ldots, E_p)\,,
$$
with $\dim E_i=n_i$, $i=1,\ldots, p$. By passing to the dual bundles,
we may assume that the functor $\phi$ is covariant in each variable.

\smallskip

Let
\begin{equation}
P(E_{\bullet})=P(E_1,E_2,\ldots,E_p)
\end{equation}
be the principal $G$-bundle associated with the bundles $E_1,E_2,\ldots,E_p$.
We define a new vector bundle:
\begin{equation}
V(E_{\bullet})=V(E_1,E_2,\ldots,E_p):=P(E_{\bullet})\times_G V\,.
\end{equation}

Suppose now that a $G$-invariant cone $\Sigma\subset V$ is given.
We set
\begin{equation}
\Sigma(E_{\bullet})=\Sigma(E_1,E_2,\ldots,E_p)
:=P(E_{\bullet})\times_{G}\Sigma\subset V(E_{\bullet})\,.
\end{equation}

We define the ``Thom polynomial'' ${\cT}^\Sigma$ to be the dual class\footnote{Compare the footnote 4 in \cite{PW}.}
of
$$
[\Sigma(R^{(1)},\ldots, R^{(p)})]\in H^*\bigl(V(R^{(1)},\ldots,
R^{(p)}),{\bf Z}\bigr)=H^*(BG,{\bf Z})\,,
$$
where $R^{(i)}$, $i=1,\dots,p$, is the pullback of the tautological vector
bundle from $BGL_{n_i}$ to
$$
BG=\prod_{i=1}^p BGL_{n_i}\,.
$$
Then, the so defined Thom polynomial
$$
{\cT}^\Sigma\in H^{*}(BG, {\bf Z})
$$
depends on the Chern classes of the $R^{(i)}$'s.

\smallskip

We shall write \ ${\cT}^{\Sigma}(E_1,\ldots,E_p)$ \ for
the Thom polynomial ${\cT}^{\Sigma}$, with $c_j(R^{(i)})$ replaced by
$c_j(E_i)$ for $i=1,\ldots,p$.

\begin{lemma}\label{LRE} For any vector bundles $E_1,E_2,\ldots,E_p$
on a variety $X$, the dual class\footnote{The meaning of the ``dual class'' for
singular $X$ is explained in \cite{PW}, Note~6.}
of $[\Sigma(E_{\bullet})]$ in
$$
H^{2\codim(\Sigma)}(V(E_{\bullet}), {\bf Z})=H^{2\codim(\Sigma)}(X,{\bf Z})
$$
is equal to ${\cT}^{\Sigma}(E_1,\ldots,E_p)$.
\end{lemma}
\proof Each $p$-tuple of bundles can be pulled back from the
universal $p$-tuple $(R^{(1)}, R^{(2)},\ldots, R^{(p)})$ of
bundles on $BG$ using a $C^{\infty}$-map. It is possible to work
entirely with the algebraic varieties and maps. One can use the
Totaro construction and representability for affine varieties
(\cite[proof of Theorem 1.3]{To}). \qed

\begin{remark}\rm
In the situation of classical Thom polynomials \cite{T}, the functor $\phi$
is the functor of $k$-jets~:
$$
(E,F)\mapsto {\cJ}^k(E,F)=\left(\bigoplus_{i=1}^k {\Sym}^iE^*\right)\otimes F\,,
$$
where $k$ is large enough, adapted to the investigated class of singularities
-- cf. \cite{PW} for details and applications.
(We note that in this situation an invariant closed subset $\Sigma$, called
in \cite{PW} a ``class of singularities'', is automatically a cone.)
\end{remark}

\section{Schur functions and the Giambelli formula}\label{Sch}

In this section, we recall the notion of {\it Schur functions}.
We also recall a geometric interpretation of them, namely the classical
{\it Giambelli formula}.

Given a partition $I=(i_1,i_2,\ldots,i_l)\in {\bf N}^l$, where
$$
i_1\ge i_2\ge \cdots\ge i_l\ge 0 \ \footnote{Since the most common references
to Schubert Calculus use weakly decreasing partitions, we follow this convention
in the present paper.}\,,
$$
and vector bundles $E$ and $F$ on a variety $X$,
the {\it Schur function}\footnote{Usually this family of functions is called
``super Schur functions'' or ``Schur functions in difference of bundles''.}
$S_I(E - F)$ is defined by the following determinant:
\begin{equation}\label{schur}
S_I(E - F)= \Bigl|
     S_{i_p-p+q}(E - F) \Bigr|_{1\leq p,q \leq l}\,,
\end{equation}
where the entries are defined by the expression
\begin{equation}\label{seg}
\sum S_i(E - F) =\prod_b(1 - b)/\prod_a(1 - a)\,.
\end{equation}
Here, the $a$'s and $b$'s are the Chern roots of $E$ and $F$
and the LHS of Eq. (\ref{seg}) is the {\it Segre class} of the virtual bundle
$E-F$. So the Schur functions $S_I(E-F)$ lie in a ring containing the Chern
classes of $E$ and $F$; e.g., we can take the cohomology ring
$H^*(X, {\bf Z})$ or the Chow ring $A^*(X)$.

Given a vector bundle $E$ and a partition $I$, we shall
write $S_I(E)$ for $S_I(E-0)$, where $0$ is the zero vector bundle.

We refer to \cite{L}, \cite{Mcd}, and \cite{PT} for the theory of Schur
functions $S_I(E)$ and $S_I(E-F)$.

\medskip

Given a smooth variety $X$, we shall identify its cohomology $H^*(X,{\bf Z})$
with its homology $H_*(X,{\bf Z})$, as is customary. More precisely, this
identification is realized via capping the cohomology classes with the
fundamental class $[X]$ of $X$, using the standard map:
$$
\cap : H^*(X,{\bf Z}) \otimes H_*(X,{\bf Z}) \to  H_*(X,{\bf Z})\,.
$$

Let $V$ be a complex vector space of dimension $N$, and let $G_m(V)$ be
the {\it Grassmannian parametrizing} $m$-dimensional subspaces of $V$. On
knows that $G_m(V)$ is a smooth projective variety of dimension $mn$,
where $n=N-m$.
We shall also use the notation $G^n(V)$ for this Grassmannian.
The Grassmannian $G_m(V)$ is stratified by Schubert cells; the closures
of these cells are Schubert varieties $\Omega_I(V_{\bullet})$,
where
$$
I=(n \ge i_1 \ge i_2 \ge \cdots \ge i_m \ge 0)
$$
is a partition, and
$$
V_{\bullet} : 0=V_0\subset V_1 \subset \cdots \subset V_N=V
$$
is a complete flag of subspaces of $V$, with $\dim {V_j}=j$
for $j=0,1,\ldots, N$.

The precise definition of $\Omega_I(V_{\bullet})$ is

\begin{equation}
\Omega_I(V_{\bullet})=\{\Lambda \in G_m(V) : \dim(\Lambda \cap V_{n+j-i_j})\ge j, \
j=1,\ldots,m \}\,.
\end{equation}
This is a subvariety of codimension $|I|=i_1+i_2+\cdots+i_m$ in $G_m(V)$.
The cohomology class $[\Omega_I(V_{\bullet})]$ does not depend
on a flag $V_{\bullet}$. We denote it by $\sigma_I$ and call
a {\it Schubert class}.

Let $Q$ denote the tautological quotient bundle on $G_m(V)$. Then
$$
\sigma_{(i)}=c_i(Q)=S_{(1,\ldots,1)}(Q)\,,
$$
where $1$ appears $i$ times, and -- more generally -- the following
{\it Giambelli formula} \cite{Giam} holds:

\begin{proposition}\label{Giamb} In the cohomology ring of $G_m(V)$, we have
\begin{equation}
\sigma_I=\Bigl|c_{i_p-p+q}(Q)\Bigr|_{1\leq p,q \leq m}=S_{I^{\sim}}(Q)\,,
\end{equation}
where
$I^{\sim}$ is the {\it conjugate} partition of $I$ (i.e. the consecutive rows of
the diagram of $I^{\sim}$ are the transposed consecutive columns of the diagram of $I$).
\end{proposition}
(Cf. \cite[Chap.~1, Sect.~5 ]{GH}, \cite[\S9.4]{F1}).

\smallskip

Given a partition $I$, consider the partition
$$
J=(n-i_m,n-i_{m-1},\ldots, n-i_1)\,.
$$
Then ({\it loc.cit.})  $\sigma_J$ is the unique Schubert class of complementary
codimension to $\sigma_I$ whose intersection with $\sigma_I$ is nonzero, and in fact
\begin{equation}\label{dual}
\int_{G_m(V)} \sigma_I \cdot \sigma_J =1\,.
\end{equation}
The class $\sigma_J$ is called the {\it complementary class} to $\sigma_I$.

\section{Cone classes for globally generated and ample vector bundles}\label{amp}

In the proof of our main result, we shall use the following
results of Fulton and Lazarsfeld from \cite{FL0}, \cite{FL}
(cf. also \cite[Chap.~12]{F}, \cite[\S8]{Laz}). Recall first some classical
definitions and
facts from \cite{F} (we shall also follow the notation from this
book). Let $E$ be a vector bundle of rank $e$ on $X$. By a {\it
cone} in $E$ we mean a subvariety of $E$ which is stable under the
natural ${\mathbb G}_{\rm m}$-action on $E$. If $C\subset E$ is a
cone of pure dimension $d$, then one may intersect its cycle $[C]$
with the zero-section of the vector bundle:
\begin{equation}
z(C,E):=s_E^*([C])\in A_{d-e}(X)\,,
\end{equation}
where $s_E^*: A_d(E)\to A_{d-e}(X)$ is the Gysin map determined by the
zero-section $s_E: X\to E$.
For a projective variety $X$, there is a well defined {\it degree} map
$\int_X: A_{0}(X)\to {\bf Z}$.

\smallskip

The following results stem from \cite[Theorem 1 (A)]{FL0}
and \cite[Theorem 2.1]{FL}.

\begin{theorem}\label{TFL}
Suppose that $E$ is a vector bundle of rank $e$ on a projective variety $X$,
and let $C\subset E$ be a cone of pure dimension $e$.

\smallskip
\noindent
(1) If some symmetric power of $E$ is globally generated, then
$$
\int_X \ z(C,E) \ge 0.
$$

\smallskip
\noindent
(2) If $E$ is ample, then
$$
\int_X \ z(C,E) > 0.
$$
\end{theorem}
Under the assumptions of the theorem,
we also have in $H_0(X,{\bf Z})$ the homology analog of $z(C,E)$,
denoted by the same symbol, and the homology degree map
$\deg_X: H_0(X,{\bf Z})\to {\bf Z}$. They are compatible with their
Chow group counterparts via the cycle map:
$A_0(X) \to H_0(X,{\bf Z})$ (cf. \cite[Chap. 19]{F}). We thus have the
same inequalities for $\deg_X \bigl(z(C,E)\bigr)$.

\section{Schur function expansions of Thom polynomials}

We follow the setting from Section \ref{Tpc}.
Since the Schur functions form an additive basis of the ring of symmetric
functions, the Thom polynomial ${\cT}^\Sigma$ is uniquely written
in the following form:
\begin{equation}\label{bul}
{\cT}^\Sigma=\sum \alpha_{I_1,\dots ,I_p} \
S_{I_1}(R^{(1)}) \ S_{I_2}(R^{(2)}) \ \cdots \ S_{I_p}(R^{(p)})\,,
\end{equation}
where $\alpha_{I_1,\dots,I_p}$ are integer coefficients.

\smallskip

We say that the functor $\phi$, associated with a representation $V$,
{\it preserves global generateness} if for a collection of globally generated
vector bundles $E_1,E_2,\dots,E_p$, the bundle
$$
\phi(E_1,E_2,\dots,E_p)
$$
is globally generated.

Examples of functors preserving global generateness over fields of
characteristic zero are {\it polynomial functors}. They are, at the same time,
quotient functors and subfunctors of the tensor power functors (cf. \cite{H}).
On the other hand, the functors: $\Hom(-,E)$ with fixed $E$, or $\Hom(-,-)$, do not
preserve global generateness.

\smallskip

The main result of the present paper is

\begin{theorem}\label{Tm} \ Suppose that the functor $\phi$
preserves global generateness. Then the coefficients $\alpha_{I_1,\ldots,I_p}$ in
Eq.~(\ref{bul}) are nonnegative. Assume additionally
that there exists a projective variety $X$ \footnote{The variety $X$ can be
singular.} of dimension greater than or equal
to $\codim(\Sigma)$, and there exist vector bundles $E_1,\ldots,E_p$ on $X$
such that the bundle $\phi(E_1,E_2,\dots,E_p)$ is ample. Then at least one
of the coefficients $\alpha_{I_1,\ldots,I_p}$ is positive.
\end{theorem}
\proof
We assume for simplicity that $p=2$ (the reasoning in general case
goes in the same way).
We want to estimate the coefficients $\alpha_{IJ}$ in the universal
expansion into products of Schur functions:
\begin{equation}
\cT^{\Sigma}(E_1,E_2)=\sum_{I,J}\alpha_{IJ} \ S_I(E_1)\cdot S_J(E_2)
\end{equation}
Let $E_1$ and $E_2$ be the pullbacks of the tautological
quotient bundles from the Grassmannians $G^{n_1}({\bf C}^{N_1})$
and $G^{n_2}({\bf C}^{N_2})$ to
$$
G^{n_1}({\bf C}^{N_1})\times G^{n_2}({\bf C}^{N_2})\,,
$$
where $N_1$ and $N_2$ are sufficiently large.
It is enough to estimate the coefficients $\alpha_{IJ}$ for such $E_1$ and $E_2$.
Let $\sigma_K \in H^*(G^{n_1}({\bf C}^{N_1}),{\bf Z})$ be the complementary class to
$\sigma_{I^{\sim}}$ and $\sigma_L \in H^*(G^{n_2}({\bf C}^{N_2}),{\bf Z})$ be the
complementary class to $\sigma_{J^{\sim}}$.
By the Giambelli formula (Proposition \ref{Giamb}) and properties of complementary
Schubert classes ({\ref{dual}), we have
$$
\alpha_{IJ}=\int_{G^{n_1}({\bf C}^{N_1})\times G^{n_2}({\bf C}^{N_2})} \cT^\Sigma(E_1,E_2)
\cdot (\sigma_K\times \sigma_L)\,.
$$
The vector bundles $E_1$ and $E_2$ are globally generated.
Hence, by the assumption, the bundle $\phi(E_1,E_2)$ is globally generated.
By Theorem \ref{TFL}(1), we thus have $\alpha_{IJ}\ge 0$.

\smallskip

Now, suppose that there exists a projective variety of dimension greater than or equal
to $\codim(\Sigma)$, and there exist vector bundles $E_1,E_2$ on $X$ such that
the bundle $\phi(E_1,E_2)$ is ample. Let $Y$ be a subvariety of $X$ of dimension
equal to $\codim(\Sigma)$.
Then, by Theorem \ref{TFL}(2), we have
$$
\int_Y\cT^{\Sigma}(E_1,E_2)>0\,.
$$
Therefore, $\cT^\Sigma\ne 0$, which implies that at least one of the coefficients
$\alpha_{IJ}$ is positive.
\qed

\smallskip

Consider now the projective variety
$$
X=\prod_{i=1}^p G^{n_i}({\bf C}^N)\,,
$$
where $N$ is sufficiently large.
We denote by $Q_i$ the pullback to $X$ of the tautological quotient bundle
on $G^{n_i}({\bf C}^N)$\,, $i=1,\ldots,p$. The bundle $Q_i$ is not ample,
but it is globally generated. Let $L$ be an ample line bundle on $X$. Then
each bundle
$$
E_i=Q_i\otimes L
$$
is ample (cf. \cite{H}).

Observe the hypotheses of the theorem are satisfied by the variety $X$,
vector bundles $E_1,\dots, E_p$, and any polynomial functor $\phi$
of positive degree. We thus obtain

\begin{corollary}\label{polfun}
If $\phi$ is a polynomial functor of positive degree, then the coefficients
$\alpha_{I_1,\dots ,I_p}$ in Eq.~(\ref{bul}) are nonnegative, and their sum
is positive.
\end{corollary}

\smallskip

In the next corollary, we use the concept of a classical Thom polynomial
associated with a map $f: M\to N$ of complex manifolds and a nontrivial
class of singularities $\Sigma$ (cf. \cite{PW}). We {\it do not}, however,
assume now that $\Sigma$ is stable.

By the theory of Schur functions, there
exist universal coefficients $\beta_{IJ}\in {\bf Z}$ such that
\begin{equation}\label{beta}
{\cT}^{\Sigma}=\sum_{I,J} \beta_{IJ} S_I(TM^*)\cdot S_J(f^*TN)\,.
\end{equation}

The following result follows from Theorem \ref{Tm}.

\begin{corollary}\label{CpS} For any pair of partitions $I,J$, we have
$\beta_{IJ}\ge 0$.
\end{corollary}

\smallskip

We also give an alternative proof of the main result from \cite{PW}.
Let $\Sigma$ be a {\it stable} singularity. Then by the Thom-Damon
theorem (\cite{T}, \cite{D}),
$$
{\cT}^{\Sigma}(c_1(M),\ldots, c_m(M),c_1(N),\ldots, c_n(N))
$$
is an universal polynomial in
$$
c_i(TM\moins f^*TN) \ \ \ \hbox{where} \ \ \
i=1,2,\ldots\,
$$
(Cf. also \cite[Theorem 2]{Ka}.)

\smallskip

Using the theory of supersymmetric functions
(cf. \cite{L}, \cite{Mcd}, \cite{PT}),
the Thom-Damon theorem can be rephrased by saying that there exist
coefficients $\alpha_I\in {\bf Z}$ such that
\begin{equation}\label{alpha}
{\cT}^{\Sigma}=\sum_I \alpha_I S_I(TM^*\moins f^*TN^*)\,,
\end{equation}
the sum is over partitions $I$ with $|I|=\codim(\Sigma)$.
The expression in Eq.~(\ref{alpha}) is unique ({\it loc.cit.}).

\begin{theorem}\label{Tpos} Let $\Sigma$ be a stable
singularity. Then for any partition $I$ the coefficient $\alpha_I$
in the Schur function expansion of the Thom polynomial
$\cT^{\Sigma}$ (cf. Eq.~(\ref{alpha})) is nonnegative.
\end{theorem}
\proof
By the theory of Schur functions ({\it loc.cit.}), we have that the coefficient of
$S_I(TM^*\moins f^*TN^*)$ in the RHS of (\ref{alpha}) is equal to the coefficient
of $S_I(TM^*)$ in the RHS of (\ref{beta}),
that is, \ $\alpha_I=\beta_{I,\emptyset}$ \ for
any partition $I$. The assertion now follows from Corollary \ref{CpS}.
\qed

\begin{remark}\rm
Note that Theorem \ref{Tm} overlaps various situations already studied
in the literature. Consider, e.g., a {\it family of quadratic forms} on
the tangent bundle of an $m$-fold $M$ with values in a line bundle $L$,
i.e. a section of
$$
{\Hom}({\Sym}^2(TM),L)\,.
$$
The singularities of such forms lead to Thom polynomials. The group
which is relevant here is $GL_m\times GL_1$ with the natural representation
in the vector space
$$
\bigoplus_{i=0}^r {\Sym}^i({\bf C}^m)\otimes
{\Hom}({\Sym}^2({\bf C}^m),{\bf C})\,.
$$
The singularity classes defined by the $0$th jet are just invariant
subsets of ${\Hom}({\Sym}^2({\bf C}^m),{\bf C})$.
The corank of the quadratic form determines the singularity class.

We recover\footnote{See also \cite{FR}.} the situation already described in
the literature in the context of {\it degeneracy loci formulas} for
morphisms with symmetries of rank $m$ bundles:
$$
E^*\to E\otimes L\,.
$$
The degrees of projective symmetric varieties were computed in
\cite{Giam1}. The Schur function formulas for a trivial $L$ were
given in \cite{HT}, \cite{JLP}. To give the formulas in full generality
\cite{P-1},
we consider, for partitions $I=(i_1,\ldots,i_m)$ and $J=(j_1\ldots,j_m)$,
the following determinant studied in the paper \cite{L0} of Lascoux:

\begin{equation}
d_{I,J}= \Bigl|{{i_a+m-a}\choose{j_b+m-b}}\Bigr|_{1\le a,b \le m}\,.
\end{equation}

Then, the Thom polynomial associated with the locus of quadratic forms
whose $\corank \ge q$ is equal to
$$
2^{-{{q}\choose {2}}} \sum_J 2^{|J|} \ d_{\rho_q,J} \ S_J(E)\cdot
S_{{{q+1}\choose{2}}-|J|}(L)\,,
$$
where $J=(j_1,\ldots,j_q)$ runs over partitions contained in the partition
$$
\rho_q=(q,q-1,\ldots,1)\,.
$$
(Cf. \cite{P-1} for details. Similar formulas are valid for antisymmetric
forms ({\it loc.cit.}).)
In particular, we obtain the positivity of the $d_{\rho_q,J}$'s --
a result known before by combinatorial methods (cf. \cite{GV}).

\smallskip

It seems to be interesting to apply Theorem \ref{Tm} to other
concrete situations, where Thom polynomials of invariant cones appear.
\end{remark}


\begin{thebibliography}{99}\small
\addcontentsline{toc}{section}{\string\numberline{}References}

\bibitem{AVGL} V. Arnold, V. Vasilev, V. Goryunov, O. Lyashko:
\emph{Singularities. Local and global theory,}
Enc. Math. Sci. Vol. {\bf 6} (Dynamical Systems VI), Springer, 1993.

\bibitem{D} J. Damon,
\emph{Thom polynomials for contact singularities,}
Ph.D. Thesis, Harvard, 1972.

\bibitem{FK} L. Feher, B. Komuves,
\emph{On second order Thom-Boardman singularities,}
Fund. Math. {\bf 191} (2006), 249--264.

\bibitem{FR} L. Feher, R. Rimanyi,
\emph{Calculation of Thom polynomials and other cohomological obstructions
for group actions,} in: ``Real and complex singularities (S\~ao Carlos
2002)'' (T. Gaffney and M. Ruas eds.), Contemporary Math. {\bf 354},
(2004), 69--93.

\bibitem{F} W. Fulton,
\emph{Intersection theory,}
Springer, 1984.

\bibitem{F1} W. Fulton,
\emph{Young tableaux,}
Cambridge University Press, 1997.

\bibitem{FL0} W. Fulton, R. Lazarsfeld,
\emph{Positivity and excess intersections,} in ``Enumerative and classical
geometry'', Nice 1981, Progress in Math. {\bf 24}, Birkha\"user (1982),
97--105.

\bibitem{FL} W. Fulton, R. Lazarsfeld,
\emph{Positive polynomials for ample vector bundles,}
Ann. of Math. {\bf 118} (1983), 35--60.

\bibitem{GV} I. Gessel, X. Viennot,
\emph{Binomial determinants, paths and hook length formulae,}
Adv. Math. {\bf 58} (1985), 300--321.

\bibitem{Giam} G. Z. Giambelli,
\emph{Risoluzione del problema degli spazi secanti,}
Mem. Accad. Sci. Torino (2) {\bf 52} (1902), 171--211.

\bibitem{Giam1} G. Z. Giambelli,
\emph{Sulle variet\'a rappresentata coll'annullare determinanti minori
contenuti in un determinante simmetrico od emisimmetrico generico di
forme,}
Atti della R. Accad. delle Scienze di Torino {\bf 44} (1906), 102--125.

\bibitem{GH} P. A. Griffiths, J. Harris,
\emph{Principles of algebraic geometry,}
Wiley \& Sons Inc., 1978.

\bibitem{HT} J. Harris, L. Tu,
\emph{On symmetric and skew-symmetric determinantal varieties,}
Topology {\bf 23} (1984), 71--84.

\bibitem{H} R. Hartshorne,
\emph{Ample vector bundles,}
Publ. Math. IHES {\bf 29} (1966), 63--94.

\bibitem{H1} R. Hartshorne,
\emph{Algebraic geometry,}
Springer, 1977.

\bibitem{JLP} T. J\'ozefiak, A. Lascoux, P. Pragacz,
\emph{Classes of determinantal varieties associated with symmetric and skew-symmetric
matrices,}
Math. USSR Izv. {\bf 18} (1982), 575--586.

\bibitem{Ka} M. E. Kazarian,
\emph{Classifying spaces of singularities and Thom polynomials,}
in: ``New developments in singularity theory'', NATO Sci. Ser. II Math.
Phys. Chem., {\bf 21}, Kluwer Acad. Publ. (2001), 117--134.

\bibitem{Kl} S. Kleiman,
\emph{The enumerative theory of singularities,}
in: ``Real and complex singularities, Oslo 1976'' (P. Holm ed.)
Sijthoff\&Noordhoff Int. Publ. (1978), 297--396.

\bibitem{L0} A. Lascoux,
\emph{Classes de Chern d'un produit tensoriel,}
C. R. Acad. Sci. Paris {\bf 286} (1978), 385--387.

\bibitem{L} A. Lascoux,
\emph{Symmetric functions and combinatorial operators on polynomials},
CBMS/AMS Lectures Notes {\bf 99}, Providence, 2003.

\bibitem{Laz} R. Lazarsfeld,
\emph{Positivity in algebraic geometry,}
Springer, 2004.

\bibitem{Mcd} I. G. Macdonald,
\emph{Symmetric functions and Hall polynomials,}
2nd edition, Oxford University Press, 1995.

\bibitem{O} O. Ozturk,
\emph{On Thom polynomials for $A_4(-)$ via Schur functions,}
Preprint, IM PAN Warszawa 2006 (670) -- to appear in Serdica Math. J.
{\bf 33} (2007).

\bibitem{P-2} P. Pragacz,
\emph{Enumerative geometry of degeneracy loci,}
Ann. Sc. Ec. Norm. Sup. {\bf 21} (1988), 413--454.

\bibitem{P-1} P. Pragacz,
\emph{Cycles of isotropic subspaces and formulas for symmetric degeneracy loci,}
in: ``Topics in algebra'' (S.~Balcerzyk et al. eds.), Banach Center Publ. {\bf 26(2)},
1990, 189--199.

\bibitem{P1} P. Pragacz,
\emph{Thom polynomials and Schur functions I,}
math.AG/0509234.

\bibitem{P2} P. Pragacz,
\emph{Thom polynomials and Schur functions: the singularities $I_{2,2}(-)$,}
Preprint MPIM Bonn 2006 (83) -- to appear in Ann. Inst. Fourier {\bf 57}
(2007).

\bibitem{P3} P. Pragacz,
\emph{Thom polynomials and Schur functions: towards the singularities
$A_i(-)$,}
Preprint MPIM Bonn 2006 (139).

\bibitem{P4} P. Pragacz,
\emph{Thom polynomials and Schur functions: the singularities $A_3(-)$,}
in preparation.

\bibitem{PT} P. Pragacz, A. Thorup,
\emph{On a Jacobi-Trudi identity for supersymmetric polynomials,}
Adv. in Math. {\bf 95} (1992), 8--17.

\bibitem{PW} P. Pragacz, A. Weber,
\emph{Positivity of Schur function expansions of Thom polynomials'}
Preprint MPIM Bonn 2006 (60), math.AG/0605308 -- to appear in Fund. Math. 
{\bf 195} (2007).

\bibitem{Rim} R. Rimanyi,
\emph{Thom polynomials, symmetries and incidences of singularities,}
Inv. Math. {\bf 143} (2001), 499--521.

\bibitem{T} R. Thom,
\emph{Les singularit\'es des applications diff\'erentiables,}
Ann. Inst. Fourier {\bf 6} (1955--56), 43--87.

\bibitem{To} B.~Totaro, \emph{The Chow ring of a classifying
space} in: ``Algebraic $K$-theory'' (W.~Raskind et al. eds.), Symp.
Pure Math. {\bf 67} (1999), AMS, 249--281.


\end{thebibliography}
\end{document}